\documentclass[12pt,leqno]{amsart} 
\newtheorem{thm}{Theorem}[section]
\newtheorem{defi}{Definition}[section]

\newtheorem{rem}{Remark}[section]
\newcommand{\be}{\begin{equation}}
\newcommand{\ee}{\end{equation}}
\newcommand{\bea}{\begin{eqnarray}}
\newcommand{\eea}{\end{eqnarray}}
\newcommand{\beb}{\begin{eqnarray*}}
\newcommand{\eeb}{\end{eqnarray*}}
\usepackage{amssymb,amsfonts,amsthm,setspace,indentfirst,mathrsfs}
\usepackage{minibox}
\numberwithin{equation}{section}
\begin{document}
%%%%%%%%%%%%%%%%%%%%%%%%%%%%%%%%%%%%%%%%%%%%%%%%%%%%%%%%%%%%%%%%%%%%%%%%%%%%%%%%%%%%%%%%%%%%%%%%%%%%%%%%%%%%%
%
\title[On curvature properties of Som-Raychaudhuri spacetime]{\bf{On curvature properties of Som-Raychaudhuri spacetime}}
\author[Absos Ali Shaikh and Haradhan Kundu]{Absos Ali Shaikh and Haradhan Kundu}
\date{\today}
\address{\noindent\newline Department of Mathematics,\newline University of 
Burdwan, Golapbag,\newline Burdwan-713104,\newline West Bengal, India}
\email{aask2003@yahoo.co.in, aashaikh@math.buruniv.ac.in}
\email{kundu.haradhan@gmail.com}
\dedicatory{}
%%%%%%%%%%%%%%%%%%%%%%%%%%%%%%%%%%%%%%%%%%%%% Abstract %%%%%%%%%%%%%%%%%%%%%%%%%%%%%%%%%%%%%%%%%%%%%%%%%%%%%
\begin{abstract}
Som-Raychaudhuri \cite{SR68} spacetime is a stationary cylindrical symmetric solution of Einstein field equation corresponding to a charged dust distribution in rigid rotation. The main object of the present paper is to investigate the curvature restricted geometric structures admitting by the Som-Raychaudhuri spacetime and it is shown that such a spacetime is a 2-quasi-Einstein, generalized Roter type,  $Ein(3)$ manifold satisfying $R.R = Q(S,R)$, $C\cdot C = \frac{2a^2}{3} Q(g,C)$, and its Ricci tensor is cyclic parallel and Riemann compatible. Finally, we make a comparison between G\"odel spacetime and Som-Raychaudhuri spacetime.
\end{abstract}
%%%%%%%%%%%%%%%%%%%%%%%%%%%%%%%%%%%%%%%%%%%%%%%%%%%%%%%%%%%%%%%%%%%%%%%%%%%%%%%%%%%%%%%%%%%%%%%%%%%%%%%%%%%%
%
\maketitle
%

%%%%%%%%%%%%%%%%%%%%%%%%%%%%%%%%%%%%%%%%%%%%%%%%%%%%%%%%%%%%%%%%%%%%%%%%%%%%%%%%%%%%%%%%%%%%%%%%%%%%%%%%%%%%%%
%																					Introduction
%%%%%%%%%%%%%%%%%%%%%%%%%%%%%%%%%%%%%%%%%%%%%%%%%%%%%%%%%%%%%%%%%%%%%%%%%%%%%%%%%%%%%%%%%%%%%%%%%%%%%%%%%%%%%%
\section{\bf Introduction}\label{intro}
%%%%%%%%%%%%%%%%%%%%%%%%%%%%%%%%%%%
Let $M$ be a connected semi-Riemannian smooth manifold of dimension $n$, $n \geq 3$, equipped with a semi-Riemannian metric $g$ of signature $(s, n-s)$, $0 \leq s \leq n$. The manifold $M$ is Riemannian if $s = 0$ or $s = n$ and is Lorentzian if $s = 1$ or $s = n-1$. Further, let $\nabla$, $R$, $S$ and $\kappa $ be respectively the Levi-Civita connection, the Riemann-Christoffel curvature tensor, the Ricci tensor and the scalar curvature of $M$.\\
\indent Physically spacetime is a mathematical model that combines space and time into a single continuum. By a spacetime we mean a four-dimensional continuum having three spatial coordinates and one temporal coordinate, in which all physical quantities are described.  In general relativity, Einstein described the gravity of a spacetime as a geometric property and linked with the geometric quantity ``curvature''. The curvature of a spacetime is directly related to its energy and momentum and the relation is specified by the famous Einstein field equation
$$
S_{\mu \nu} - \frac{\kappa}{2}g_{\mu \nu} + g_{\mu \nu} \Lambda = \frac{8 \pi G}{c^4} T_{\mu \nu},
$$
where $\Lambda$ is the cosmological constant, $G$  is Newton's gravitational constant, $c$ is the speed of light in vacuum and $T_{\mu \nu}$ is the energy-momentum tensor. Hence to investigate the nature of a spacetime one have to study its geometry, especially, the curvature of the spacetime and the curvature restricted geometric structures admitting by the spacetime.\\
%=============================================================================
\indent Symmetry plays the key role to describe the geometry of a manifold. In 1926 Cartan defined the notion of local symmetry \cite{Cart26} by $\nabla R = 0$, but many spacetimes are not locally symmetric and hence it is necessary to generalize such notion. During the last eight decades many researchers are working to generalize the notion of local symmetry and describe the geometry of corresponding spacetimes. In consequence, there arose various generalizations of local symmetry, such as semisymmetry (\cite{Cart46}, \cite{Szab82}), pseudosymmetry (in the sense of Deszcz \cite{AD83} and Chaki \cite{Chak87}), weak symmetry \cite{TB89}, recurrent manifold \cite{Walk50} etc.\\
%=========================================================================
\indent In 1968 Som and Raychaudhuri \cite{SR68} presented a family of stationary cylindrical symmetric solution of Einstein field equation, known as Som-Raychaudhuri solution, corresponding to a charged dust distribution in rigid rotation. In terms of cylindrical coordinates $(t, r, z, \phi)$, the stationary cylindrical symmetric metric is given by
\be\label{eq1.1}
ds^2 = g_{00} dt^2 - e^{2\psi} (dr^2+dz^2)-l d\phi^2 + 2m d\phi dt,
\ee
where $g_{00}, \psi, l$ and $m$ are functions of $r$ alone. According to Som and Raychaudhuri's consideration $g_{00}=1$, $m=a r^2$, $l=r^2-a^2r^4$, $\psi =0$, $a = constant (\neq 0)$, and hence the Som-Raychaudhuri metric is given by
\be\label{eq1.2}
ds^2 = dt^2 - (r^2-a^2r^4) d\phi^2 - dr^2 - dz^2 + 2 a r^2 d\phi dt.
\ee
Thus Som-Raychaudhuri metric can be written as a G\"odel type metric (see \cite{DHJKS14}, \cite{RT83} and references therein) as follows:
\be\label{eq1.3}
ds^{2} = (dt + h(r)\, d\phi )^{2} - (f(r))^{2}\, d\phi ^{2} - dr^{2} - dz^{2},
\ee
where $h(r) = a r^{2}$, $a = constant (\neq 0)$ and $f(r) = r$. Again for $h(r) \, =\, \frac{2 \sqrt{2}}{m} \, sinh^{2} (\frac{m r}{2})$ and $f(r)\, =  \, \frac{2}{m} \, sinh (\frac{m r}{2})\, cosh (\frac{m r}{2})$ ($m=constant$), the metric \eqref{eq1.3} reduces to the G\"{o}del metric \cite{Gode49} (see \cite{RT83}, eq. (1.6)). In terms of Cartesian coordinates $(x,y,z,t)$, the G\"odel metric can be presented as
\be\label{eq1.4}
ds^2 = (dt)^2 + \frac{1}{2}e^{2 m x}(dy)^2 - (dx)^2 - (dz)^2 + 2e^{m x}dt dy.
\ee
%============================================================================
\indent The main object of the present paper is to investigate the curvature restricted geometric structures admitting by the Som-Raychaudhuri spacetime. The paper is organized as follows. Section 2 deals with the definitions of various curvature restricted geometric structures. In Section 3 we investigate the curvature restricted geometric structures admitting by Som-Raychaudhuri spacetime. Among others, it is shown that Som-Raychaudhuri spacetime is $2$-quasi-Einstein, Ricci generalized pseudosymmetric and is of generalized Roter type. Section 4 is concerned with the comparison between Som-Raychaudhuri spacetime and G\"odel spacetime.
%%%%%%%%%%%%%%%%%%%%%%%%%%%%%%%%%%%%%%%%%%%%%%%%%%%%%%%%%%%%%%%%%%%%%%%%%%%%%%%%%%%%%%%%%%%%%%%%%%%%%%%%%%%%%%%%%%%%
%                                               Geometric structures
%%%%%%%%%%%%%%%%%%%%%%%%%%%%%%%%%%%%%%%%%%%%%%%%%%%%%%%%%%%%%%%%%%%%%%%%%%%%%%%%%%%%%%%%%%%%%%%%%%%%%%%%%%%%%%%%%%%%
\section{\bf Curvature restricted geometric structures}
%%%%%%%%%%%%%%%%%%%%%%%%%%%%%%%%%%%%%%%%%%%%%%%%%%%%%%%
In this section we define various curvature restricted geometric structures. Let us consider a connected semi-Riemannian smooth manifold $M$ of dimension $n$, $n \geq 3$, equipped with a semi-Riemannian metric $g$. Let $C^{\infty}(M)$, $\chi(M)$ and $\chi^*(M)$ be the algebra of all smooth functions, the Lie algebra of all smooth vector fields and the Lie algebra of all smooth 1-forms on $M$ respectively. We define $C^{\infty}(M)$-linear endomorphisms over $\chi(M)$, namely, $X\wedge_A Y$, $\mathcal{R}(X,Y)$, $\mathcal{C}(X,Y)$, $\mathcal{P}(X,Y)$, $\mathcal{W}(X,Y)$ and $\mathcal {K}(X,Y)$ as follows (\cite{DHJKS14}, \cite{SKgrt}):
\begin{eqnarray*}
(X\wedge_A Y)X_1 &=& A(Y,X_1)X-A(X,X_1)Y,\\
\mathcal{R}(X,Y)X_1 &=& [\nabla_X,\nabla_Y]X_1-\nabla_{[X,Y]}X_1,\\
\mathcal{C}(X,Y) &=& \mathcal{R}(X,Y)-\frac{1}{n-2}(X\wedge_g \mathcal{S} Y + \mathcal{S} X\wedge_g Y - \frac{\kappa }{n-1}X\wedge_g Y),\\
\mathcal{P}(X,Y) &=& \mathcal{R}(X,Y)-\frac{1}{n-1}X\wedge_S Y,\\
\mathcal{W}(X,Y) &=& \mathcal{R}(X,Y) - \frac{\kappa }{n (n-1)} X\wedge_g Y \ \ \mbox{and}\\
\mathcal{K}(X,Y) &=& {\mathcal{R}}(X,Y) - \frac{1}{n-2}(X \wedge _{g} {\mathcal{S}}Y + {\mathcal{S}}X \wedge _{g} Y ),
\end{eqnarray*}
where $X,Y,X_1\in\chi(M)$, $A$ is a symmetric $(0,2)$-tensor and $\mathcal{S}$ is the Ricci operator, defined
by $g(X,\mathcal{S} Y) = S(X,Y)$. Throughout the paper we consider $X, Y, X_i, Y_i \in \chi(M)$, $i = 1,2, \ldots$. In view of the above endomorphisms we have various curvature tensors given as follows (\cite{DHJKS14}, \cite{SK14}, \cite{GLOG2},\cite{Ishii}):
\begin{eqnarray*}
&&\mbox{Gaussian curvature tensor}: G(X_1,X_2,X_3,X_4) = g((X_1\wedge_g X_2)X_3,X_4),\\
&&\mbox{Riemann-Christoffel curvature tensor}: R(X_1,X_2,X_3,X_4) = g(\mathcal{R}(X_1,X_2)X_3,X_4),\\
&&\mbox{conformal curvature tensor}: C(X_1,X_2,X_3,X_4) = g(\mathcal{C}(X_1,X_2)X_3,X_4),\\
&&\mbox{projective curvature tensor}: P(X_1,X_2,X_3,X_4) = g(\mathcal{P}(X_1,X_2)X_3,X_4),\\
&&\mbox{concircular curvature tensor}: W(X_1,X_2,X_3,X_4) = g(\mathcal{W}(X_1,X_2)X_3,X_4),\\
&&\mbox{conharmonic curvature tensor}: K(X_1,X_2,X_3,X_4) = g(\mathcal{K}(X_1,X_2)X_3,X_4).
\end{eqnarray*}
\indent Again there are three $(0,2)$-tensors, namely, $S^2$, $S^3$ and $S^4$, defined by
$$S^{2}(X,Y) = S(X,\mathcal{S} Y), \ S^{3}(X,Y) = S(X,\mathcal{S}\mathcal{S} Y) \ \mbox{and} \ S^{4}(X,Y) = S(X,\mathcal{S}\mathcal{S}\mathcal{S} Y)$$
respectively, and called Ricci tensor of level 2, 3 and 4 respectively.\\
%=============================================
\indent One can easily operate an $C^{\infty}(M)$-linear endomorphism $\mathcal H$ on a $(0,k)$-tensor $T$, $k\geq 1$, and get the tensor $\mathcal H T$ given by (\cite{DG02}, \cite{DGHS98}, \cite{DH03})
\beb\label{hdot}
(\mathcal{H} T)(X_1,X_2,\cdots,X_k) = -T(\mathcal{H}X_1,X_2,\cdots,X_k) - \cdots - T(X_1,X_2,\cdots,\mathcal{H}X_k).
\eeb
In particular, for $\mathcal H = \mathcal R(X,Y)$ we get a $(0,k+2)$-tensor $R\cdot T$ (see \cite{SDHJK15}, \cite{SK14} and also references therein) as follows:
$$R\cdot T(X_1,X_2,\cdots,X_k,X,Y) = -T(\mathcal R(X,Y)X_1,X_2,\cdots,X_k) - \cdots - T(X_1,X_2,\cdots,\mathcal R(X,Y)X_k).$$
Similarly for $\mathcal H = \mathcal C(X,Y)$ (resp., $\mathcal W(X,Y)$ and $\mathcal K(X,Y)$) we get the $(0,k+2)$-tensor $C\cdot T$ (resp., $W\cdot T$ and $K\cdot T$). If $A$ is a symmetric (0,2)-tensor, then for $\mathcal H = X\wedge_A Y$, we get a $(0,k+2)$-tensor $Q(A,T)$ (\cite{SDHJK15}, \cite{SK14}, \cite{Tach74}) given by
\beb
&&Q(A,T)(X_1,X_2, \ldots ,X_k,X,Y) = ((X \wedge_A Y)\cdot T)(X_1,X_2, \ldots ,X_k)\\
&&= A(X, X_1) T(Y,X_2,\cdots,X_k) + \cdots + A(X, X_k) T(X_1,X_2,\cdots,Y)\\
&& - A(Y, X_1) T(X,X_2,\cdots,X_k) - \cdots - A(Y, X_k) T(X_1,X_2,\cdots,X).
\eeb
%===============================================
\begin{defi}$($\cite{AD83}, \cite{Cart46}, \cite{Desz92}, \cite{SK14}, \cite{Szab82}$)$ 
A semi-Riemannian manifold $M$ is said to be $T$-semisymmetric type if $\mathcal H T = 0$ and it is said to be $T$-pseudosymmetric type if $\left(\sum\limits_{i=1}^k c_i \mathcal H_i\right) T = 0$, where $\mathcal H, \mathcal H_i$, $i=1,\ldots, k$, $(k\ge 2)$, are $C^{\infty}(M)$-linear endomorphisms and $c_i\in C^{\infty}(M)$, called the associated scalars.
\end{defi}
%=======
\indent In particular, if $\mathcal H = \mathcal R(X,Y)$ and $T=R$ (resp., $S$, $C$, $W$ and $K$), then $M$ is called semisymmetric (resp., Ricci, conformally, concircularly and conharmonically semisymmetric). Again, if $i =2$, $\mathcal H_1 = \mathcal R(X,Y)$, $\mathcal H_2 = X\wedge_g Y$ and $T= R$ (resp., $S$, $C$, $W$ and $K$), then $M$ is called Deszcz pseudosymmetric (resp., Ricci, conformally, concircularly and conharmonically pseudosymmetric). Especially, if $i =2$, $\mathcal H_1 = \mathcal C(X,Y)$, $\mathcal H_2 = X\wedge_g Y$ and $T =C$, then $M$ is called a manifold of pseudosymmetric Weyl conformal curvature tensor. Again, if $i =2$, $\mathcal H_1 = \mathcal R(X,Y)$, $\mathcal H_2 = X\wedge_S Y$ and $T =R$, then $M$ is called Ricci generalized pseudosymmetric \cite{DD91}.
%===============================================
\begin{defi}
A semi-Riemannian manifold is said to be a $k$-quasi-Einstein manifold if $rank (S - \alpha g) = k$, $0\le k\le (n-1)$, for a scalar $\alpha$. In particular, for $k=0$, $S$ is linearly dependent with $g$ and the manifold is Einstein. Again, for $k=1$, the manifold is called quasi-Einstein and in this case we have $S = \alpha g + \beta \Pi \otimes \Pi$, $\Pi\in\chi^*(M)$ and $\alpha, \beta \in C^{\infty}(M)$, $\otimes$ denotes the tensor product. If $\alpha = 0$, then a quasi-Einstein manifold is called Ricci simple.
\end{defi}
\indent We note that every Einstein manifold is quasi-Einstein and every $(k-1)$-quasi-Einstein manifold is $k$-quasi-Einstein for $1\le k\le (n-1)$. To obtain the explicit form of Ricci tensor, there are various generalizations of quasi-Einstein manifold.
%====================================
\begin{defi}
A semi-Riemannian manifold $M$ is said to be generalized quasi-Einstein in the sense of Chaki \cite{CHAKI} (resp., in the sense of De and Ghosh \cite{UG1}) if
$$S = \alpha g + \beta \Pi \otimes \Pi + \gamma [\Pi \otimes \Phi + \Phi \otimes \Pi]$$
$$(\mbox{resp., } \ S = \alpha g + \beta \Pi \otimes \Pi + \gamma \Phi \otimes \Phi)$$
holds on $U = \{x\in M : (S-\nu_1 g-\nu_2 \xi\otimes \xi)_x \neq 0 \mbox{ for any scalars $\nu_1, \nu_2$ and any 1-form $\xi$}\}$ for some $\alpha, \beta, \gamma \in C^{\infty}(M)$ and $\Pi, \Phi \in \chi^*(M)$ such that their associated vector fields are mutually orthogonal.
\end{defi}
%============================
\indent We note that the notion of generalized quasi-Einstein manifolds in the sense of Chaki is different from the that by De and Ghosh.
%============================
 Again in 2009, Shaikh \cite{AS1} introduced the notion of pseudo quasi-Einstein manifold and studied its geometric properties and relativistic significance along with the existence of such notion by several non-trivial examples.
\begin{defi} \cite{AS1}
A semi-Riemannian manifold $M$ is said to be pseudo quasi-Einstein if
\be\label{pqe}
S = \alpha g + \beta \Pi \otimes \Pi + \gamma E
\ee
holds on $U$ for some $\alpha, \beta, \gamma \in C^{\infty}(M)$, $\Pi\in \chi^*(M)$ and a trace free $(0,2)$-tensor $E$ such that the associated vector field $V$ corresponding to $\Pi$ satisfies $E(X,V) =0$.
\end{defi}
\indent In particular, if $E = \Phi \otimes \Phi$ (resp., $\Pi \otimes \Phi + \Phi \otimes \Pi$) then a pseudo quasi-Einstein manifold reduces to a generalized quasi-Einstein manifold by De and Ghosh (resp., by Chaki).\\
%==============================================
\indent For two $(0,2)$-tensors $A$ and $E$, their Kulkarni-Nomizu product $A\wedge E$ is given by (\cite{DGHS11}, \cite{GLOG2}, \cite{Glog02}, \cite{SK14})
\begin{eqnarray*}
(A\wedge E)(X_1,X_2,X,Y)&=&A(X_1,Y)E(X_2,X)+A(X_2,X)E(X_1,Y)\\
&-&A(X_1,X)E(X_2,Y)-A(X_2,Y)E(X_1,X).
\end{eqnarray*}
\begin{defi} 
A semi-Riemannian manifold $M$ is said to be a Roter type manifold $($\cite{Desz03}, \cite{Desz03a}$)$ $($resp., generalized Roter type manifold $($\cite{SDHJK15}, \cite{SKgrt}, \cite{SKgrtw}$))$ if
$$R = N_1 g\wedge g + N_2 g\wedge S + N_3 S\wedge S \ \mbox{ and}$$
$$(\mbox{resp., } R = L_1 g\wedge g + L_2 g\wedge S + L_3 S\wedge S + L_4 g\wedge S^2 + L_5 S\wedge S^2 + L_6 S^2\wedge S^2)$$
holds for some $N_i, L_j \in C^{\infty}(M)$, $1\le i\le 3$ and $1\le j\le 6$.
\end{defi}
%============================
\begin{defi}(\cite{Bess87}, \cite{SKgrt})
A semi-Riemannian manifold $M$ is said to be $Ein(2)$, $Ein(3)$ and $Ein(4)$ if
\be
S^2 + a_1 S + a_2 g = 0,
\ee
\be
S^3 + a_3 S^2 + a_4 S + a_5 g = 0 \ \mbox{and}
\ee
\be
S^4 + a_6 S^3 + a_7 S^2 + a_8 S + a_9 g = 0
\ee
holds respectively, for some $a_i \in C^{\infty}(M)$, $1\le i \le 9$.
\end{defi}
%============================================
\indent In 1978 Gray presented two new classes of manifolds, namely, manifolds of Codazzi type Ricci tensor and manifolds of cyclic parallel Ricci tensor, lies between the class of Ricci symmetric manifolds and the class of manifolds of constant scalar curvature.
\begin{defi}
A semi-Riemannian manifold $M$ is said to be of Codazzi type Ricci tensor (resp. cyclic parallel Ricci tensor) (see, \cite{DHJKS14}, \cite{GRAY} and references therein), if
\[(\nabla_{X_1} S)(X_2, X_3) = (\nabla_{X_2} S)(X_1, X_3)\]
\[\left(\mbox{resp.} \ \ \ (\nabla_{X_1} S)(X_2, X_3) + (\nabla_{X_2} S)(X_3, X_1) + (\nabla_{X_3} S)(X_1, X_2) = 0 \ \right)\]
holds on $M$.
\end{defi}
%============================================
\begin{defi}\label{def2.8}
A symmetric $(0, 2)$-tensor $E$ on a semi-Riemannian $M$ is said to be Riemann compatible or simply $R$-compatible (\cite{MM12}, \cite{MM12a}) if
\[
R(\mathcal E X_1, X,X_2,X_3) + R(\mathcal E X_2, X,X_3,X_1) + R(\mathcal E X_3, X,X_1,X_2) = 0,
\]
holds, where $\mathcal E$ is the endomorphism corresponding to $E$ defined as $g(\mathcal E X_1, X_2) = E(X_1, X_2)$. Again a 1-form $\Pi$ is said to be Riemann compatible if $\Pi\otimes \Pi$ is Riemann compatible.
\end{defi}
Similarly we can define conformal compatibility (also known as Weyl compatibility see, \cite{DGJPZ13} and \cite{MM13}), concircular compatibility and conharmonic compatibility.
%%%%%%%%%%%%%%%%%%%%%%%%%%%%%%%%%%%%%%%%%%%%%%%%%%%%%%%%%%%%%%%%%%%%%%%%%%%%%%%%%%%%%%%%%%%%%%%%%%%%%%%%%%%%%%%%%%%%%%%%%
%                                                  Main Results
%%%%%%%%%%%%%%%%%%%%%%%%%%%%%%%%%%%%%%%%%%%%%%%%%%%%%%%%%%%%%%%%%%%%%%%%%%%%%%%%%%%%%%%%%%%%%%%%%%%%%%%%%%%%%%%%%%%%%%%%%
\section{\bf Som-Raychaudhuri spacetime admitting geometric structures}\label{main}
%%%%%%%%%%%%%%%%%%%%%%%%%%%%%%%%%%%%%%%%%%%%%%%%%%%%%%%%%%%%%%%%%%%%%%%%%%%%%%%%%%%
Let us consider $\mathbb R^4$ equipped with the Som-Raychaudhuri metric $g$ which, in terms of cylindrical coordinates $(t,\phi,r,z)$, is given by
\be\label{met}
ds^{2} = g_{ij}dx^i dx^j = (dt +  a r^{2}\, d\phi )^{2}- r^{2}\, d\phi ^{2}- dr^{2}- dz^{2},
\ee
where the coordinates $t,\phi,r$ and $z$ are respectively designated as $x^1,x^2,x^3$ and $x^4$.\\
%===========================================================
\indent From \eqref{met} the non-zero components (upto symmetry) of the $R$, $S$, $\nabla R$ and $\nabla S$ are given by
$$R_{1212}=-a^2 r^2, \ \ R_{1313}=-a^2, \ \ R_{1323}= -a^3 r^2, \ \ R_{2323}=-a^2 r^2 \left(a^2 r^2+3\right);$$
$$S_{11}= S_{33}=-2 a^2, \ \ S_{12}=-2 a^3 r^2, \ \ S_{22}=-2 \left(a^4 r^4+a^2 r^2\right);$$
$$R_{1223,2}= -4 a^3 r^3, \ \ R_{1323,3}= -4 a^3 r, \ \ R_{2323,3}=-8 a^4 r^3;$$
$$S_{12,3}=- S_{13,2}= -4 a^3 r, \ \ -\frac{1}{2}S_{22,3}= S_{23,2}=4 a^4 r^3.$$
%
%=============================================================
\indent The value of the scalar curvature of the Som-Raychaudhuri metric is $2 a^2$ and the non-zero components (upto symmetry) of $C$, $W$ and $K$ are given by
$$C_{1212}=-\frac{2}{3} a^2 r^2, \ \ -C_{1313}= \frac{1}{2}C_{1414}= C_{3434}=\frac{2 a^2}{3}, \ \ C_{2323}=-\frac{2}{3} a^2 r^2 \left(a^2 r^2+2\right),$$
$$-C_{1323}= \frac{1}{2}C_{1424}= \frac{2}{3} a^3 r^2, \ \ C_{2424}=\frac{2}{3} \left(2 a^4 r^4+a^2 r^2\right);$$
%-------------
$$W_{1212}=-\frac{7}{6} a^2 r^2, \ \ -\frac{1}{7}W_{1313}= -W_{1414}= W_{3434}=\frac{a^2}{6}, \ \ 
W_{2323}=-\frac{1}{6} a^2 r^2 \left(7 a^2 r^2+17\right),$$
$$\frac{1}{7}W_{1323}= W_{1424}= -\frac{1}{6} a^3 r^2, \ \ 
W_{2424}=-\frac{1}{6} a^2 r^2 \left(a^2 r^2-1\right);$$
%-----------
$$K_{1212}=-a^2 r^2, \ \ -K_{1313}= K_{1414}= K_{3434}=a^2,$$
$$-K_{1323}= K_{1424}= a^3 r^2, \ \ K_{2323}=-K_{2424}=-a^2 r^2 \left(a^2 r^2+1\right).$$
%
%==================================================================
\indent Again the non-zero components (upto symmetry) of 
$R\cdot R$, $R\cdot S$, $R \cdot C$, 
$C \cdot R$, $C \cdot C$, 
$Q(g, R)$, $Q(S, R)$, $Q(g, C)$ and $Q(S, C)$ are given by
%---------------------------------------------------------------------------------------------------------
$$-R\cdot R_{122313}= R\cdot R_{132312}=4 a^4 r^2, \ \ -2R\cdot R_{122323}= R\cdot R_{232312}=8 a^5 r^4;$$
%----------------------------------------------------------------------------------------------------------
$$R \cdot S_{1212}=4 a^4 r^2, \ \ R \cdot S_{1313}=4 a^4, \ \ R \cdot S_{2212}=8 a^5 r^4,$$
$$R \cdot S_{1323}= R \cdot S_{2313}=4 a^5 r^2, \ \ R \cdot S_{2323}=4 a^6 r^4;$$
%----------------------------------------------------------------------------------------------------------
$$-R\cdot C_{122313}= R\cdot C_{132312}= -R\cdot C_{142412}=2 a^4 r^2,$$
$$-2R\cdot C_{122323}= R\cdot C_{232312}= -R\cdot C_{242412}=4 a^5 r^4, \ \ -R\cdot C_{143413}=2 a^4,$$
$$-R\cdot C_{143423}= -R\cdot C_{243413}=2 a^5 r^2, \ \ -R\cdot C_{243423}=2 a^6 r^4;$$
%----------------------------------------------------------------------------------------------------------
$$4 C \cdot R_{121424}= -C \cdot R_{122313}= 2 C \cdot R_{122414}= C \cdot R_{132312}=\frac{8 a^4 r^2}{3},$$
$$-2 C \cdot R_{122323}= \frac{3}{8} C \cdot R_{122424}= C \cdot R_{232312}=\frac{16 a^5 r^4}{3},$$
$$2 C \cdot R_{131434}= C \cdot R_{133414}=\frac{4 a^4}{3}, \ \ C \cdot R_{232434}=\frac{2}{3} a^4 r^2 \left(a^2 r^2+3\right),$$
$$2 C \cdot R_{132434}= C \cdot R_{133424}= 2 C \cdot R_{142334}= C \cdot R_{233414}=\frac{4 a^5 r^2}{3}, \ \ C \cdot R_{233424}=\frac{4 a^6 r^4}{3}-2 a^4 r^2;$$
%----------------------------------------------------------------------------------------------------------
$$\frac{3}{8}C\cdot S_{1212}=a^4 r^2, \ \ \frac{3}{8}C\cdot S_{1313}= -\frac{3}{8}C\cdot S_{1414}= \frac{3}{4}C\cdot S_{3434}=a^4,$$
$$\frac{3}{8}C\cdot S_{1323}= -\frac{3}{8}C\cdot S_{1424}=a^5 r^2, \ \ \frac{3}{8}C\cdot S_{2323}=a^6 r^4, \ \ -\frac{3}{4}C\cdot S_{2424}=a^4 r^2 \left(2 a^2 r^2-1\right);$$
%----------------------------------------------------------------------------------------------------------
$$\frac{3}{4}C \cdot C_{121424}= -\frac{3}{4}C \cdot C_{122313}= \frac{3}{4}C \cdot C_{132312}= -\frac{3}{4}C \cdot C_{142412}= -\frac{3}{4}C \cdot C_{233424}=a^4 r^2,$$
$$-\frac{3}{4}C \cdot C_{122323}= \frac{3}{4}C \cdot C_{122424}= \frac{3}{8}C \cdot C_{232312}= -\frac{3}{8}C \cdot C_{242412}=a^5 r^4,$$
$$\frac{3}{4}C \cdot C_{132434}= \frac{3}{4}C \cdot C_{142334}= -\frac{3}{4}C \cdot C_{143423}= -\frac{3}{4}C \cdot C_{243413}=a^5 r^2,$$
$$\frac{3}{4}C \cdot C_{131434}= -\frac{3}{4}C \cdot C_{143413}=a^4, \ \ \frac{3}{4}C \cdot C_{232434}=a^4 r^2 \left(a^2 r^2+1\right), \ \ -\frac{3}{4}C \cdot C_{243423}=a^6 r^4;$$
%----------------------------------------------------------------------------------------------------------
%----------------------------------------------------------------------------------------------------------
$$-4Q(g,R)_{121424}= Q(g,R)_{122313}= 4Q(g,R)_{122414}= -Q(g,R)_{132312}=4 a^2 r^2,$$
$$2Q(g,R)_{122323}= -Q(g,R)_{232312}=8 a^3 r^4, \ \ -Q(g,R)_{131434}= Q(g,R)_{133414}=a^2,$$
$$-Q(g,R)_{132434}= Q(g,R)_{133424}= -Q(g,R)_{142334}= Q(g,R)_{233414}=a^3 r^2,$$
$$-Q(g,R)_{232434}= Q(g,R)_{233424}=a^2 r^2 \left(a^2 r^2+3\right);$$
%----------------------------------------------------------------------------------------------------------
$$-Q(S,R)_{122313}= Q(S,R)_{132312}=4 a^4 r^2, \ \ -2Q(S,R)_{122323}= Q(S,R)_{232312}=8 a^5 r^4;$$
%----------------------------------------------------------------------------------------------------------
$$-Q(g,C)_{121424}= Q(g,C)_{122313}= -Q(g,C)_{132312}= Q(g,C)_{142412}= Q(g,C)_{233424}=2 a^2 r^2,$$
$$2Q(g,C)_{122323}= -2Q(g,C)_{122424}= -Q(g,C)_{232312}= Q(g,C)_{242412}=4 a^3 r^4,$$
$$-Q(g,C)_{232434}=2 a^2 r^2 \left(a^2 r^2+1\right), \ \ Q(g,C)_{243423}=2 a^4 r^4, \ \ -Q(g,C)_{131434}= Q(g,C)_{143413}=2 a^2,$$
$$-Q(g,C)_{132434}= -Q(g,C)_{142334}= Q(g,C)_{143423}= Q(g,C)_{243413}=2 a^3 r^2;$$
%----------------------------------------------------------------------------------------------------------
$$-\frac{1}{2}Q(S,C)_{121424}= -Q(S,C)_{122313}= Q(S,C)_{122414}= Q(S,C)_{132312}= Q(S,C)_{142412}=\frac{4 a^4 r^2}{3},$$
$$-2Q(S,C)_{122323}= -2Q(S,C)_{122424}= Q(S,C)_{232312}= Q(S,C)_{242412}=\frac{8 a^5 r^4}{3},$$
$$-\frac{1}{2}Q(S,C)_{131434}= Q(S,C)_{133414}= Q(S,C)_{143413}=\frac{4 a^4}{3}, \ \ -Q(S,C)_{232434}=\frac{4}{3} a^4 r^2 \left(2 a^2 r^2+1\right)$$
$$-\frac{1}{2}Q(S,C)_{132434}= Q(S,C)_{133424}= -\frac{1}{2}Q(S,C)_{142334}= Q(S,C)_{143423}= Q(S,C)_{233414}$$
$$= Q(S,C)_{243413}=\frac{4 a^5 r^2}{3}, \ \ \ Q(S,C)_{233424}=\frac{4}{3} a^4 r^2 \left(a^2 r^2+1\right), \ \ Q(S,C)_{243423}=\frac{4 a^6 r^4}{3}.$$
%----------------------------------------------------------------------------------------------------------
\indent Hence from the above, it follows that Som-Raychaudhuri spacetime has the following curvature properties:\\
(i) Not quasi-Einstein but 2-quasi-Einstein as
$$S = 2a^2 g + (\Pi\otimes \Phi + \Phi\otimes \Pi),$$
where $\Pi = (1,ar^2,0,1/\sqrt{2})$, $\Phi = (-2a^2,-2a^3r^2,0,\sqrt{2}a^2)$, and thus $rank(S-2a^2 g) = 2$.\\
(ii) Ricci tensor is not of Codazzi type but cyclic parallel.\\
(iii) Generalized quasi-Einstein manifold in the sense of Chaki for
$$\alpha = 2a^2, \ \ \beta= -12a^4r^4, \ \ \gamma=1, \ \ \Pi =  \left(\frac{1}{ar^2},1,0,-\frac{1}{\sqrt{2}ar^2}\right), \ \ \Psi = \left(4a^3r^2,4a^4r^4,0,4\sqrt{2}a^3r^2\right).$$
(iv) Generalized quasi-Einstein manifold in the sense of De and Ghosh for
$$\alpha = 2a^2, \ \ \beta= 1, \ \ \gamma=-1, \ \ \Pi =  \left(0,0,0,\sqrt{2}a\right), \ \ \Psi = \left(2a,2a^2r^2,0,0\right).$$
(v) Pseudo quasi-Einstein manifold for
$$\alpha = \frac{2}{3}a^2, \ \ \beta= 1, \ \ \gamma=-1, \ \ \Pi =  \left(0,0,0,\sqrt{\frac{2}{3}}a\right) \ \ \mbox{and}$$
$$E= \left(
\begin{array}{cccc}
 \frac{8}{3}a^2 & \frac{8}{3}a^3r^2 & 0 & 0 \\
 \frac{8}{3}a^3r^2 & \frac{4}{3}a^2r^2(1+2a^2r^2) & 0 & 0 \\
 0 & 0 & \frac{4}{3}a^2 & 0 \\
 0 & 0 & 0 & 0
\end{array}
\right).$$
(vi) Not pseudosymmetric but special Ricci generalized pseudosymmetric ($R\cdot R = Q(S,R)$).\\
(vii) Not conformally pseudosymmetric but of pseudosymmetric Weyl conformal curvature tensor satisfying $C\cdot C = \frac{2a^2}{3} Q(g,C)$.\\
(viii) Not Roter type but satisfies the generalized Roter type condition
$$R = L_1 g\wedge g + L_3 S\wedge S -\frac{L_1}{2 a^4} g\wedge S^2 + \left(\frac{1}{4 a^4}-\frac{L_3}{a^2}\right) S\wedge S^2
 + \left(\frac{L_1}{16 a^8}-\frac{1}{32 a^6}+\frac{L_3}{4 a^4}\right) S^2\wedge S^2,$$
where $L_{1}, L_{3}$ are arbitrary scalars.\\
(ix) Not $Ein(2)$ but an $Ein(3)$ manifold since $S^3 = 4a^4 S$.\\
(x) Also satisfies the pseudosymmetric type condition
$$L_{11} R\cdot R + L_{13} R\cdot C + L_{13} C\cdot R + L_{14} C\cdot C = L_{11} Q(S,R) + \left(-\frac{5 a^2 L_{13}}{3}-\frac{2 a^2 L_{14}}{3}\right) Q(g,C) + L_{13} Q(S,C),$$
where $L_{11}, L_{13}, L_{14}$ are arbitrary scalars.\\
(xi) The Ricci tensor is Riemann compatible as well as conformal compatible, concircular compatible and conharmonic compatible.\\
%========================================================
Hence we can state the following:
\begin{thm}
The Som-Raychaudhuri spacetime $(M,g)$ is (i) of cyclic parallel Ricci tensor, (ii) 2-quasi-Einstein, (iii) $Ein(3)$ (iv) generalized quasi-Einstein both in the sense of Chaki, and of De and Ghosh, (v) pseudo quasi-Einstein (vi) special Ricci generalized pseudosymmetric, (vii) manifold of pseudosymmetric Weyl conformal curvature, (viii) generalized Roter type and (ix) its Ricci tensor is Riemann compatible, conformal compatible, concircular compatible and conharmonic compatible.
\end{thm}
%========================================================
%%%%%%%%%%%%%%%%%%%%%%%%%%%%%%%%%%%%%%%%%%%%%%%%%%%%%%%%%%%%%%%%%%%%%%%%%%%%%%%%%%%%%%%%%%%%%%%%%%%%%%%%%%%%%%%%%%%%%%%%%
%                                                  Extension
%%%%%%%%%%%%%%%%%%%%%%%%%%%%%%%%%%%%%%%%%%%%%%%%%%%%%%%%%%%%%%%%%%%%%%%%%%%%%%%%%%%%%%%%%%%%%%%%%%%%%%%%%%%%%%%%%%%%%%%%%
\section{\bf Som-Raychaudhuri spacetime as a G\"odel type spacetime}\label{extension}
%%%%%%%%%%%%%%%%%%%%%%%%%%%%%%%%%%%%%%%%%%%%%%%%%%%%%%%%%%%%%%%%%%%%%%%%%%%%%%%%%%%%%
Recently Deszcz et al. \cite{DHJKS14} studied the geometric properties of G\"odel metric and showed that the G\"odel spacetime is Ricci simple. We recall that both the G\"odel metric and Som-Raychaudhuri metric are of G\"odel type. In this section, after presenting some geometric properties of G\"odel type metric, we make a comparison between Som-Raychaudhuri spacetime and G\"odel spacetime.\\
%-----------------------------------------------------------
\indent Let $M$ be a connected smooth manifold endowed with the G\"odel type metric
\be\label{gtm}
ds^{2} = (dt + h(r)\, d\phi )^{2} - (f(r))^{2}\, d\phi ^{2} - dr^{2} - dz^{2},
\ee
in cylindrical coordinates $(t, r, z, \phi)$. Then the non-zero components of $R$, $S$ and $C$ are given by
%-----------------------------------------------
$$R_{1212}=-\frac{\left(h'\right)^2}{4}, \ \ R_{1313}=-\frac{\left(h'\right)^2}{4 f^2}, \ \ R_{1323}=-\frac{2 f^2 h''-2 f f' h'+h \left(h'\right)^2}{4 f^2},$$
$$R_{2323}=\frac{-4 f^2 h h''-3 f^2 \left(h'\right)^2+4 f h f' h'+4 f^3 f''-h^2 \left(h'\right)^2}{4 f^2},$$
$$S_{11}=-\frac{\left(h'\right)^2}{2 f^2}, \ \ S_{12}=-\frac{f^2 h''-f f' h'+h \left(h'\right)^2}{2 f^2},$$
$$S_{22}=\frac{-2 f^2 h h''-f^2 \left(h'\right)^2+2 f h f' h'+2 f^3 f''-h^2 \left(h'\right)^2}{2 f^2}, \ \ S_{33}=\frac{2 f f''-\left(h'\right)^2}{2 f^2},$$
$$C_{1212}=-\frac{1}{6} \left(\left(h'\right)^2-f f''\right), \ \ C_{1313}= -\frac{1}{2}C_{1414}= -C_{3434}=\frac{f f''-\left(h'\right)^2}{6 f^2},$$

$$C_{1323}=-\frac{3 f^2 h''-2 f h f''-3 f f' h'+2 h \left(h'\right)^2}{12 f^2}, \ \ C_{1424}=\frac{3 f^2 h''-4 f h f''-3 f f' h'+4 h \left(h'\right)^2}{12 f^2},$$

$$C_{2323}=\frac{-3 f^2 h h''+f h^2 f''-2 f^2 \left(h'\right)^2+3 f h f' h'+2 f^3 f''-h^2 \left(h'\right)^2}{6 f^2},$$
$$C_{2424}=-\frac{-3 f^2 h h''+2 f h^2 f''-f^2 \left(h'\right)^2+3 f h f' h'+f^3 f''-2 h^2 \left(h'\right)^2}{6 f^2},$$
where $f' =\frac{df}{d r}$, $f'' =  \frac{d f'}{d r}$, $h' =\frac{dh}{d r}$, $h'' =  \frac{d h'}{d r}$, $h^{(3)} =  \frac{d h''}{d r}$ and $h^{(4)} =  \frac{d h^{(3)}}{d r}$.\\
%---------------------------------------------
From above we can state the following:
\begin{thm}
Let $M$ be a connected smooth manifold endowed with the G\"odel type metric \eqref{gtm}. Then we have the following:\\
(i) $M$ is a 3-quasi-Einstein manifold and it is a 2-quasi Einstein manifold if $f h'' = f' h'$,  and Ricci simple manifold if $h'^2 = 2 f f''$.\\
(ii) $M$ is a special Ricci generalized pseudosymmetric manifold (i.e., $R\cdot R = Q(S,R)$), and it is a semisymmetric manifold if $f=c h'$ and $h' = c^2 h'''$ for a constant $c$.\\
(iii) Its Ricci tensor is Riemann compatible, Weyl compatible, concircular compatible and conharmonic compatible, and the Ricci tensor is cyclic parallel if $h^{(4)} h'-h^{(3)} h''=0$ and $f=c h'$ for a constant $c$.\\
(iv) $M$ satisfies $C\cdot C= L Q(g,C)$ if $f =c h'$ for a constant $c$, and in this case $L = \frac{f f''-h'^2}{6 f^2} = -\frac{\left(h'-a^2 h^{(3)}\right)}{6 a^2 h'}$.\\
(v) If $\tau = ( h'^{2} - 2 f f'' ) ( f^{2} h''^{2} - 2 f f' h' h'' - h'^{4} + 2 f f'' h'^{2} + f'^{2} h'^{2})$  is non-zero at every point of $M$, then $M$ satisfies the generalized Roter type condition \cite{DHJKS14}
$$R = L_{1} S \wedge S + L_{2} S \wedge S^{2} + L_{3} S^{2} \wedge S^{2},$$
for 
$L_{1} = \frac{f^{2}}{\tau } (2 f^{2} h''^{2} - 4 f f' h' h'' - 3 h'^{4} + 8 f f'' h'^{2} + 2 f'^{2} h'^{2} - 8 f^{2} f''^{2}), L_2 =\frac{2 f^{4}}{\tau} ( h'^{2} - 4 f f'')$ and $L_{3} = - \frac{4 f^{6}}{\tau }$ . 
\end{thm}
%%%%%%%%%%%%%%%%%%%%%%%%%%%%%%%%%%%%%%%%%%%%%%%%%%%%%%%%%%%%%%%%%%%%%%%%%%%%%%%%%%%%%%%%%%%%%%%%%%%%%%%%%%%%%%%%%%%%%%
We have now state some additional information about G\"odel metric.
\begin{rem}
Let $M$ be a connected smooth manifold endowed with the G\"odel metric \eqref{eq1.4}. Then we have the following:\\
(i) The 1-form $\omega = (0, m e^{m r},0,m)$ is Riemann compatible as well as Weyl compatible, concircular compatible and conharmonic compatible.\\
(ii) The Ricci tensor is also Riemann compatible, Weyl compatible, concircular compatible and conharmonic compatible.\\
(iii) $M$ is an Ein(2) manifold satisfying $S^2= m^2 S$.\\
(iv) It satisfies the pseudosymmetric type condition $P\cdot R = \frac{2}{3} Q(S,R)$.
\end{rem}
%%%%%%%%%%%%%%%%%%%%%%%%%%%%%%%%%%%%%%%%%%%%%%%%%%%%%%%%%%%%%%%%%%%%%%%%%%%%%%%%%%%%%%%%%%%%%%%%%%%%%%%%%%%%%%%%%%%%%%
Hence we have the following comparison between G\"odel spacetime and Som-Raychaudhuri spacetime.\\
\textbf{A. Similarity:}\\
(i) Their Ricci tensors are cyclic parallel but not of Codazzi type.\\
(ii) They are of special Ricci generalized pseudosymmetric.\\
(iii) They are of pseudosymmetric Weyl conformal curvature tensor.\\
(iv) Their Ricci tensors are Riemann compatible as well as conformal compatible, concircular compatible and conharmonic compatible.\\
\textbf{B. Dissimilarity:}\\
(i) G\"odel spacetime is quasi-Einstein but Som-Raychaudhuri spacetime is proper 2-quasi-Einstein.\\
(ii) G\"odel spacetime satisfies $K.K = 0$ but Som-Raychaudhuri spacetime satisfies $K\cdot K = -a^2 Q(g,K)$.\\
(iii) G\"odel spacetime is not generalized Roter type (since in this case $\tau=0$) but Som-Raychaudhuri spacetime is generalized Roter type.\\
(iv) G\"odel spacetime is Ein(2) but Som-Raychaudhuri spacetime is Ein(3).
%%%%%%%%%%%%%%%%%%%%%%%%%%%%%%%%%%%%%%%%%%%%%%%%%%%%%%%%%%%%%%%%%%%%%%%%%%%
\begin{thm}
Both the G\"odel spacetime and Som-Raychaudhuri spacetime are Ricci generalized pseudosymmetric, pseudosymmetric Weyl conformal curvature tensor, and their Ricci tensors are Riemann, conformal, concircular and conharmonic compatible, and cyclic parallel but not of Codazzi type.
\end{thm}
\begin{thm}
G\"odel spacetime is a quasi-Einstein, Ein(2) manifold satisfies $K.K = 0$ but Som-Raychaudhuri spacetime is a proper 2-quasi-Einstein, Ein(3) manifold satisfies $K\cdot K = -a^2 Q(g,K)$. Moreover Som-Raychaudhuri spacetime is generalized Roter type but G\"odel spacetime is not generalized Roter type.
\end{thm}
%%%%%%%%%%%%%%%%%%%%%%%%%%%%%%%%%%%%%%%%%%%%%%%%%%%%%%%%%%%%%%%%%%%%%%%%%%%%%%%%%%%%%%%%%%%%%%%%%%%%%%%%%%%%%%
%                                                  Conclusion
%%%%%%%%%%%%%%%%%%%%%%%%%%%%%%%%%%%%%%%%%%%%%%%%%%%%%%%%%%%%%%%%%%%%%%%%%%%%%%%%%%%%%%%%%%%%%%%%%%%%%%%%%%%%%%%%
\section{\bf Conclusion}
%=========================
Som-Raychaudhuri spacetime is a G\"odel type stationary cylindrical symmetric solution of Einstein field equation corresponding to a charged dust distribution in rigid rotation. In the present paper we investigate the curvature restricted geometric structures admitting by the Som-Raychaudhuri spacetime and it is shown that such a spacetime is a 2-quasi-Einstein, pseudo quasi-Einstein, generalized quasi-Einstein, generalized Roter type, $Ein(3)$ manifold and satisfies $R.R = Q(S,R)$, $C\cdot C = \frac{2a^2}{3} Q(g,C)$, and its Ricci tensor is cyclic parallel and Riemann compatible. We obtain sufficient conditions for a G\"odel type metric to be 2-quasi Einstein, cyclic Ricci parallel and semisymmetric. The comparisons between the G\"odel spacetime and Som-Raychaudhuri spacetime are also discussed. It is shown that both the G\"odel spacetime and Som-Raychaudhuri spacetime are Ricci generalized pseudosymmetric, pseudosymmetric Weyl conformal curvature tensor, and their Ricci tensors are Riemann compatible and cyclic parallel. It is also shown that G\"odel spacetime is a quasi-Einstein, not generalized Roter type, Ein(2) manifold and satisfies $K.K = 0$, but Som-Raychaudhuri spacetime is a proper 2-quasi-Einstein, generalized Roter type, Ein(3) manifold and fulfills $K\cdot K = -a^2 Q(g,K)$.\\
%
%===============================================================================================================
\noindent\newline
\textbf{Acknowledgment:} 
The second named author gratefully acknowledges to CSIR, New Delhi (File No. 09/025 (0194)/2010-EMR-I) for the financial assistance. All the algebraic computations of Section \ref{main} are performed by a program in Wolfram Mathematica.
%%%%%%%%%%%%%%%%%%%%%%%%%%%%%%%%%%%%%%%%%%%%%%%%%%%%%%%%%%%%%%%%%%%%%%%%%%%%%%%%%%%%%%%%%%%%%%%%%%%%%%%%%%%%%%%%%
%%%%%%%%%%%%%%%%%%%%%%%%%%%%%%%%%%%%%%%%%%%%%%%%%%%%%%%%%%%%%%%%%%%%%%%%%%%%%%%%%%%%%%%%%%%%%%%%%%%%%%%%%%%%%%%%%

%%%%%%%%%%%%%%%%%%%%%

%%%%%%%%%%%%%%%%%%%%%%%%%%%%%%%%%%%%%%%%%%%%%%%%%%%%%%%%%%%%%%%%%%%%%%%%%%%%%%%%%%%%%%%%%%%%%%%%%%%%
\end{document}